\documentclass[12pt]{amsart}
\usepackage{amscd,amssymb,times,amssymb}
\begin{document}

\title{Geometric versus homotopy theoretic equivariant bordism}
\author{Bernhard Hanke}

\begin{abstract}
By results of L\"offler and Comeza\~na, the Pontrjagin-Thom map 
from geometric $G$-equivariant bordism to homotopy theoretic 
equivariant bordism is  
injective for compact abelian $G$. If $G = S^1 \times \ldots \times S^1$,  
we prove  that the  associated fixed point square 
is a pull back square,  thus   
confirming a recent conjecture of Sinha \cite{S2}. 
This is used in order to determine the image of the Pontrjagin-Thom 
map for toral $G$. 
\end{abstract}

\maketitle

\newtheorem{thm}{Theorem}
\newtheorem{prop}{Proposition}
\newtheorem{lem}{Lemma}
\newtheorem{cor}{Corollary}
\newtheorem{conj}[thm]{Conjecture}

\theoremstyle{remark}
\newtheorem{ex}{Example}
\newtheorem{rem}{Remark}

\theoremstyle{definition}
\newtheorem{defn}{Definition}

\newcommand{\G}{\Gamma}
\newcommand{\blowup}{\overline{\bC P^2}}

\newcommand{\Spi}{Spin}
\newcommand{\SO}{SO}
\newcommand{\SU}{SU}
\newcommand{\MU}{MU}
\newcommand{\BU}{BU}
\newcommand{\BO}{BO}
\newcommand{\EU}{EU}
\newcommand{\Spc}{Spin^c}
\newcommand{\GL}{GL}
\newcommand{\MSO}{MSO}
\newcommand{\MSpc}{MSpin^c}
\newcommand{\MSpi}{MSpin}
\newcommand{\K}{K} 
\newcommand{\RO}{RO}
\newcommand{\muu}{mu}

\newcommand{\im}{{\rm im}}
\newcommand{\N}{\mathbb{N}}
\newcommand{\C}{\mathbb{ C}}
\newcommand{\F}{\mathbb{ F}}
\newcommand{\bfq}{{\mathbb {Q}}}
\newcommand{\bfz}{{\mathbb {Z}}}
\newcommand{\bfn}{{\mathbb {N}}}
\newcommand{\Z}{\mathbb{ Z}}
\newcommand{\del}{\partial}
\newcommand{\Hom}{\operatorname{Hom}}
\newcommand{\hra}{\hookrightarrow}
\newcommand{\lra}{\longrightarrow}
\newcommand{\MM}{\Gamma}
\newcommand{\PM}{\mathbb{ M}}
\newcommand{\GG}{\mathcal{ G}}
\newcommand{\II}{\mathcal{ I}}
\newcommand{\BB}{\mathcal{ B}}
\newcommand{\OO}{\mathcal{ O}}
\newcommand{\Sl}{\mathcal{ S}}
\newcommand{\Q}{\mathbb{ Q}}
\newcommand{\R}{\mathbb{ R}}
\newcommand{\HH}{\mathbb{ H}}
\newcommand{\rk}{\operatorname{rk}}
\newcommand{\RP}{\mathbb{ R}\mathbb{ P}}
\newcommand{\ra}{{\rightarrow}}
\newcommand{\Sym}{\operatorname{Sym}}
\newcommand{\cp}{{\C P^2}}
\newcommand{\first}{\cp\#\blowup}
\newcommand{\second}{\cp\#2\blowup}
\newcommand{\quadric}{S^2 \x S^2}
\newcommand{\x}{\times}
\newcommand{\CC}{\mathcal{ C}}
\newcommand{\Ha}{\mathcal{ H}}

\section{Introduction} Geometric equivariant bordism theory was introduced  
by Conner and Floyd 
\cite{CF} as an important tool for studying groups of transformations 
on smooth manifolds. Later, the homotopy theoretic analogue was described
by tom Dieck \cite{tD1} and in a slightly different way by 
Br\"ocker and Hook \cite{BrHo} in order to develop a conceptual
approach to  the  localization 
theorem  of Atiyah-Segal \cite{AtSe} and to study equivariant characteristic
numbers of $G$-manifolds \cite{tD1a,tD2}. The homotopy theoretic definition 
 turned out to be the correct description of 
equivariant bordism in the context of equivariant stable homotopy theory,
because it is stable under suspensions with arbitrary 
$G$-representations. It has gained recent interest in the 
context of general completion
theorems \cite{GM} analogous to the
 Atiyah-Segal completion theorem in 
equivariant $K$-theory and in connection
with the study of equivariant formal group laws \cite{CGK,CGK2}.

The coefficients of equivariant bordism theories are difficult to calculate 
and only partial results are known.
Generators for the geometric equivariant 
bordism rings have been found in \cite {KY} for $G = S^1$ and 
in \cite{K} for $G = \Z / p $. The semifree geometric 
$S^1$-bordism ring is calculated in \cite{S2}. In \cite{Kriz}, the homotopy theoretic 
$\Z/p$-bordism ring is described in terms of a pull back square involving 
localization with respect to nontrivial Euler classes. The  
paper \cite{S1} investigates generators and relations 
in the homotopy theoretic $S^1$-bordism ring.

Let $G$ be a compact Lie group. We denote by   
$\Omega_*^{G}$ the geometric unitary $G$-bordism ring and by $\MU^{G}_*$ the
homotopy theoretic unitary $G$-bordism ring (for 
definitions see Section \ref{review} 
below).  Both theories are 
regarded as graded over the integers.

Similar to the nonequivariant case there is a Pontrjagin-Thom map
\[
    \Psi :  \Omega_*^{G}  \to \MU^{G}_* \, . 
\]
However, the proof of 
$\Psi$ being an isomorphism fails if $G \neq \{ e\}$ due to the 
lack of generic equivariant transversality \cite{Pe}.

If $G = \Z/p$ is of prime order, tom Dieck 
\cite{tD1} showed that $\Psi$ is injective. This result 
was generalized by L\"offler \cite{L1} and Comeza\~na \cite{C} 
to the case that $G$ is a compact abelian Lie group. It is not known, if 
$\Psi$ is injective in general.

On the other hand, if $G$ is nontrivial, then 
$\MU^G_*$ contains nonzero elements in negative  degrees: 
Each complex $G$-representation $W$ gives rise to an Euler class
\[
  e_W \in \MU^{G}_{-2|W|}
\]
which is zero if and only if $W^G \neq 0$. Here and in the following $|W|$ 
denotes the complex dimension of a complex $G$-representation $W$. 
In particular, $\Psi$ is not surjective, if $G \neq \{e\}$, since  
the geometrically defined theory $\Omega^G_*$ is concentrated in
nonnegative degrees. 

Our paper contributes to the understanding of the relation 
between geometric and homotopy theoretic equivariant bordism
if $G$ is a torus of the form $S^1 \times \ldots \times S^1$. On the one 
hand it complements the classical 
work of tom Dieck, L\"offler and others. 
On the other hand it should provide some useful information 
in connection with the present interest in equivariant stable 
homotopy theory.

 In the recent paper \cite{S2} the image of $\Psi$ is described for semifree 
$S^1$-bordism and $\Z / p $-bordism. 
Conjecture 4.1.~in {\em loc.~cit.}~gives a conjectural description 
of $\im(\Psi)$ for $G = S^1$.

For toral $G$, in  Theorem \ref{main} we show the existence of a pull back square 
\[
   \begin{CD} 
           \Omega_*^G @>>> \MU_*[e_V^{-1}, Y_{V,d}]     \\
            @V \Psi VV        @V {\rm incl.} VV \\
           \MU_*^G @>>> \MU_*[e_V, e_V^{-1}, Y_{V,d}]  \, .   
   \end{CD}
\]
It refines \cite{tD1}, Proposition 4.1, and generalizes \cite{S2}, 
Corollary 2.11. In particular, it confirms the 
aforementioned conjecture. The horizontal maps in this diagram capture the 
normal data around fixed point 
sets of elements in $\Omega_*^G$ and $\MU^G_*$. These 
normal data are expressed in terms of polynomials with $\MU_*$-coefficients
in Euler classes $e_V$ of irreducible complex $G$-representations ($V$ 
running through a complete set of isomorphism classes of 
nontrivial irreducible complex $G$-representations), their formal
inverses and in certain classes $Y_{V,d}$, $2 \leq d < \infty$, of 
degree $2d$. We can think of $e_V^{-1}$ as being 
geometrically represented by the $G$-bundle 
\[
   V \to *
\]
and - more generally - of 
the classes $Y_{V,d}$ as being represented by the $G$-bundles 
\[
   E \otimes V \to \C P^{d-1}
\]
where $E \to \C P^{d-1}$ is the hyperplane line bundle (i.e. the 
normal bundle of $\C P^{d-1}$ in $\C P^{d}$).

The horizontal maps in the above 
diagram are injective (assuming that $G$ is a torus) and we can conclude
that the  Pontrjagin-Thom map induces a ring isomorphism 
\[
   \Omega^G_* \cong \MU^G_* \cap \MU_*[e_V^{-1}, Y_{V,d}] \, , 
\]
this intersection taking place in $\MU_*[e_V, e_V^{-1}, Y_{V,d}]$.  
In particular, the realizability of an element in $\MU_*^G$ as an actual 
$G$-manifold is 
expressed in terms of a representation theoretic condition on the  normal 
data around the fixed point set. This should be contrasted with 
the usual failure of  $G$-equivariant transversality as discussed in 
\cite{Pe}: In general $G$-equivariant transversality is 
obstructed not only by  local invariants expressible in terms 
of Euler classes, but also by global ones, see Theorem 
1 in {\em loc.~cit.} In 
the context of  equivariant bordism however,
the vanishing of the local obstructions turns out to be  sufficient for 
equivariant transversality within the equivariant 
homotopy class of a given stable 
map $S^W \to T(\xi^G_n)$ with respect to the zero section of $T(\xi^G_n)$. 
Here, $T(\xi^G_n)$ 
denotes a certain equivariant Thom space (see below) and $S^W$ is 
the one point compactification $W \cup \{\infty\}$ with base point $\infty$.

If $G$ is not a torus and not of prime order, geometric
realizability can not be read off from normal data
around the fixed point set. We will illustrate this
by easy examples. 

In the case of toral $G$, we 
refine the description of geometric $G$-bordism using 
the fact that the image of the fixed point map 
\[
   \MU^G_* \to \MU_*[e_V, e_V^{-1} ,  Y_{V,d}]
\]
is characterized by  
integrality conditions \cite{tD1}. Details will be given 
below. 

The methods used in the proof of Theorem \ref{main} 
are mainly classical, based on considering bordism 
with respect to families of subgroups and equivariant characteristic 
numbers.

\section{Review of equivariant bordism} \label{review}

We briefly recall the definitions and constructions that 
are necessary in order to formulate our 
results. If not stated otherwise, $G$ denotes 
a compact Lie group. 

\begin{defn} \label{geo} Let $M$ be a closed smooth $G$-manifold. 
We define a {\em stable almost complex $G$-structure} on $M$ to be a 
complex $G$-structure on 
\[
    TM \oplus \underline{\R}^k \to M 
\]
for some $k$. Here, the trivial bundle 
\[  
   \underline{\R}^k = M \times \R^k \to M
\]
is equipped with the trivial $G$-action on the fibres $\R^k$ and 
by a complex $G$-structure, we mean a $G$-equivariant  
bundle map
\[
  J :   TM \oplus \underline{\R}^k \to TM \oplus \underline{\R}^k 
\]  
over $M$ which satisfies $J^2 = -1$. Two  
stable almost complex $G$-structures are identified, if after stabilization 
with further $\underline{\C}$-summands (equipped with trivial $G$-actions
on the fibres), the induced complex $G$-structures are $G$-homotopic through 
complex $G$-structures. 
\end{defn}

Note that with this definition, the (unstabilized) normal 
bundle of the $G$-embedding 
\[
    M^G \subset M
\]
of the fixed point set has an induced structure of a complex $G$-bundle. 
This property is fundamental for many calculations in geometric 
equivariant bordism theory. 

We spelt out Definition \ref{geo} in order to avoid 
possible confusion with other variants of stable almost complex 
$G$-structures. For example, \cite{C} contains
the notion of a normally almost 
complex $G$-structure on a smooth $G$-manifold $M$, being 
defined as a complex $G$-bundle structure on
the stable normal bundle of the embedding of $M$ in some real 
$G$-representation, see Definition 2.1 in {\it loc.~cit}.
The following example illustrates the difference between this 
notion and ours.  

Let $S^2$ be equipped with the antipodal $\Z/2$-action. 
Consider the $\Z/2$-embedding 
\[
  S^2 \subset \R^3
\]
where $\Z/2$ acts on $\R^3$ by multiplication with $-1$. 
The equivariant normal bundle of this embedding is 
\[   
    \nu = S^2 \times \R \to S^2
\]
with the trivial action on $\R$. Hence, $S^2$ is 
a normally almost complex $\Z/2$-manifold in the sense of the 
above definition,  
although the quotient $S^2/ (\Z/2)$ is not even orientable.
It is easy to check that $S^2$ equipped with this $\Z/2$-action does 
not admit a stable almost complex $\Z/2$-structure in the sense of 
Definition \ref{geo}. Furthermore, we remark that for a free 
stable almost complex  $G$-manifold $M$, the quotient $M/G$ always has
 an induced stable almost complex structure.  
From now on, we will exclusively work with Definition \ref{geo}.   

For a $G$-space $X$, we define the geometric unitary $G$-equivariant 
bordism groups of $X$, denoted $\Omega^G_n(X)$,
in the usual way as $G$-bordism  classes of singular stable almost complex $G$-manifolds $M^n \to X$. The coefficients 
\[
   \Omega^G_* := \Omega^G_*({\rm pt.})
\]
of this theory are equipped with a ring structure induced
by the diagonal $G$-action on the cartesian product of two $G$-manifolds. 

We recall the definition of homotopy theoretic
unitary $G$-bordism introduced in \cite{tD1}. Modern 
expositions of equivariant stable homotopy theory and 
$\RO(G)$-graded homology theories can be 
found in \cite{Ca,tD4,LMS}. For the sake of brevity, we restrict 
ourselves to  an ad hoc definition of the relevant structures. Let
\[
    \xi^G_n \to \BU(n,G)
\]
be the universal unitary $n$-dimensional Grassmannian $G$-bundle  and $T(\xi^G_{n})$ its Thom space regarded as a pointed $G$-space. 

Let $X$ be a pointed $G$-CW complex. By definition, 
\[  
 \widetilde{\MU}^G_{2n}(X) := \lim_{\longrightarrow_W} [S^W, T(\xi^G_{|W|-n}) \wedge X]^G
\]
where $[-,-]^G$ is the group of homotopy classes of pointed 
$G$-equivariant maps
and where the colimit is with respect to a  directed set of unitary
$G$-representations which ultimately contains each irreducible 
with arbitrary multiplicity. Recall (\cite{tD4}, Lemma (II.6.1)) 
that the suspension maps in the induced directed system of 
homotopy groups are uniquely determined, because we are using complex representations. Furthermore, we set  
\[
   \widetilde{MU}^G_{2n-1}(X) := \widetilde{\MU}^G_{2n}(S^1 \wedge X)
\]
where $S^1$ (with basepoint $1$) carries the trivial $G$-action. Again, the 
coefficients 
\[
    \MU^G_* : = \widetilde{MU}^G_{*}(S^0)
\]
carry a canonical ring structure. 

The Pontrjagin-Thom map 
\[
   \Psi: \Omega^G_* \to \MU^G_* 
\]
is defined as follows: Let a stable almost complex $G$-manifold $M$ represent a 
bordism class in $\Omega^G_{2n}$.  We embed  $M$ in a 
unitary $G$-representation $W$,  and let  $\nu$  
be the normal bundle of this embedding. If $W$ is chosen large enough, 
the bundle $\nu$ carries an induced structure of a unitary $G$-bundle 
and we construct a map  
\[
     S^W \to T(\xi^G_{|\nu|}) 
\]
from the universal bundle map  $\nu \to \xi^G_{|\nu|}$  by collapsing the complement of the disk bundle 
$D(\nu)$ to the base 
point in $T(\xi^G_{|\nu|})$. We explain in some detail how the required 
unitary $G$-bundle structure on $\nu$ is obtained. Let
\[
    \iota : M \to V 
\]
be a $G$-embedding into a unitary $G$-representation $V$ with 
(real) normal $G$-bundle $\mu$. We
have fixed the structure of a unitary $G$-bundle on 
\[
   TM \oplus \underline{\R}^{2k}
\]
for some $k$ (recall that we assume $\dim M$ to be even). Let $N$ be 
(the stably unique)  complementary unitary $G$-bundle, i.e. 
\[ 
    (TM \oplus \underline{\R}^{2k}) \oplus N = M \times U
\]
with a certain unitary $G$-representation $U$. The bundle 
$\mu \oplus \underline{U}$ is the normal bundle of the embedding 
\[
    M \to V \times U \, , ~ x \mapsto (\iota(x), 0)
\]
and has an induced unitary $G$-structure, because
\[
 \mu \oplus \underline{U} =  \mu \oplus TM \oplus \underline{\C}^k \oplus N
 = \underline{V} \oplus \underline{\C}^k \oplus N \, . 
\]

If $M$ is odd-dimensional, we  obtain a map 
\[  
      S^{W} \to S^1 \wedge T(\xi^G_{|\widetilde{\nu}|}) 
\]
because $\nu$ has a direct sum decomposition 
$\nu = \underline{\R} \oplus \widetilde{\nu}$ with a unitary $G$-bundle $\widetilde{\nu}$
(after choosing $W$ large enough).  

The Pontrjagin-Thom map is a map of graded $\MU_*$-algebras. We cite 

\begin{prop}[\cite{C}, Theorem 5.4] \label{loffler} For 
a compact abelian Lie group $G$, 
the Pontrjagin-Thom map is a split monomorphism of $\MU_*$-modules. 
\end{prop} 

A quick proof of the injectivity of $\Psi$ for the case $G = S^1 \times 
\ldots \times S^1$ can be given, if one uses the fact 
(see \cite{tD2}, Theorem 3) that 
for topologically cyclic $G$, the characteristic 
number map 
\[
    \Omega^G_{*} \to K_G^{-*}[[a_1, a_2, \ldots ]]\, , 
\]
which factors through the Pontrjagin-Thom map, is injective.

Let $W$ be a complex $G$-representation of dimension $n$ and let 
\[
   \begin{CD} 
       W  @>>> \xi^G_{n} \\
       @VVV          @VVV \\
        * @>>> \BU(n,G)
   \end{CD}
\]
be the universal bundle map. The induced
map $S^0 \subset S^W \to T(\xi^G_n)$ represents the Euler class  
\[
   e_W \in \MU_{-2n}^G
\]
which is zero if and only if $W$ contains a trivial $G$-representation
as a direct summand. In particular, for nontrivial $G$, the Pontrjagin-Thom 
map is not surjective, since the geometric theory $\Omega^{G}_*$ 
is concentrated in nonnegative degrees. 

For dealing with the different $G$-representations we need a certain 
amount of bookkeeping device. We follow the 
exposition in \cite{tD1}.  Let $J$ be a set containing exactly one representative of 
each isomorphism class of nontrivial irreducible complex 
$G$-representations. We define the graded ring
\[
    A_*(G) := \Z [ \Z J ] 
\]
where $\Z J$ is the free abelian group generated by $J$. We 
consider $\Z J$ as a subgroup of the abelian group underlying the 
complex representation ring $R(G)$.  
The grading of $A_*(G)$ is induced by the real dimension 
of the underlying (virtual) $G$-representations. There is  an isomorphism 
of graded rings
\[
   A_*(G) \cong \Z[e_V, e_{V}^{-1}] \, , 
\]
where $V$ runs over elements in $J$ and each $e_V$ has 
degree $-2|V|$. This isomorphism is induced by the map  
\[
   \Z J \ni \sum_{V \in J} \alpha_V V \mapsto \prod e_V^{-\alpha_V} \, . 
\]
Elements in the ring $A_*(G)$ 
will be used in order to describe  the 
(virtual) $G$-representations occuring around 
the fixed set of a given element in $\MU^G_*$. 

Fix a basepoint $1$ in $\BU$ and let 
\[
   B \subset \BU^J
\]
be the subset consisting of $J$-indexed families in $\BU$ 
of which only finitely many components are different 
from the base point. Whitney sum of vector bundles induces an 
$H$-space structure 
\[
    m : \BU \times \BU \to \BU \, .
\]
We can assume that $m(1,1)=1$. 
With the componentwise multiplication, we 
get an $H$-space structure 
on $B$. The space $B$ is the classifying space for unitary 
$G$-bundles over base spaces $X$ with trivial $G$-action and with $G$
acting without 
fixed points in each fibre.  More precisely, the 
$V$-th component of a map 
\[
   X \to B
\]
classifies the $V$-isotypical 
component $E_V$ of a given such $G$-bundle 
\[
     \bigoplus_{V \in J} E_V \otimes V \to X \, .
\]
In \cite{tD1}, p. 350, tom Dieck 
constructs a map of graded $\MU_*$-algebras
\[
    \phi_{\MU}:= \phi : \MU^G_* \to \MU_*(B) \otimes A_*(G) \, .
\]
Heuristically, this map is given by restriction to fixed point 
sets. This becomes especially clear in the slightly different description 
of the map $\phi$ in \cite{S1}, which we now briefly 
recall (for a more detailed exposition, see \cite{S1}): Let  
\[
    f : S^W \to T(\xi^G_n)
\]
represent an element $ c\in \MU^G_{2(|W|-n)}$. Consider the 
restriction 
\[
    f^G : (S^W)^G \to T(\xi^G_n)^G 
\]
of $f$ to fixed point sets. This map represents an element in 
\[
   (\phi^G \MU^G)_{2(|W| - n)} \, , 
\]
where $\phi^G \MU^G$ is the geometric fixed points spectrum associated with
the equivariant unitary bordism spectrum $\MU^G$. We have an equivalence 
of ring spectra 
\[
  \phi^G \MU^G \simeq {\bf I}_{R(G)} \wedge \MU \wedge B_+ 
\]
where $\MU$ is the usual unitary bordism spectrum and 
\[
    {\bf I}_{R(G)} : = \bigvee_{W \in R(G), ~ |W| = 0} S^{2|W^G|} 
\]
is the one point union of suspended sphere spectra equipped 
with a ring spectrum 
structure induced by addition of elements in $R(G)$, cf. Theorem 4.9. in \cite{S1}. 
We consider the stable homotopy of ${\bf I}_{R(G)} \wedge \MU \wedge B_+$
as a collection of suspended copies of $\MU_*(B)$ and hence get an 
isomorphism of graded rings 
\[
  \omega :  ({\bf I}_{R(G)} \wedge \MU \wedge B_+)_* \cong \MU_*(B) \otimes A_*(G) \, . 
\]
It is induced by identifying the summand of $\MU_*(B)$ indexed by
$W \in R(G)$ with  
\[
    \MU_*(B) \otimes \big( e_{(U^G)^{\perp}} \cdot (e_{(V^G)^{\perp}}^{-1})
   \big) \subset \MU_*(B) \otimes A_*(G) \, .
\]
Here, we write $W = U - V$ with unitary $G$-representations $U$ and $V$ and 
use splittings $U = U^G \oplus (U^G)^{\perp}$ and 
$V = V^G \oplus (V^G)^{\perp}$. 
Finally, we set
\[
    \phi(c) := \omega([f^G]) \, .
\]
For elements in $\MU^G_*$ of odd degree, the definition of $\phi$
is similar. Note the tautological equation
\[
    \phi_{\MU}(e_V) = e_V 
\]
for each nontrivial irreducible $G$-representation $V$. 

The K\"unneth formula shows that for each finite subset $J' \subset J$, we have
\[
   \MU_*(\prod_{V \in J'} \BU) \cong \bigotimes_{V \in J'} \MU_*(\BU)\, . 
\]
where the tensor product is over $\MU_*$. By a standard application of the 
Atiyah-Hirzebruch spectral sequence,   
\[
   \MU_*(\BU) \cong \MU_*[X_1, X_2, X_3 , \ldots ]
\]
is a polynomial algebra over $\MU_*$ in infinitely many indeterminates $X_d$ of 
degree $2d$, $ 1 \leq d< \infty$, see \cite{Kochman}. We can and will 
choose the polynomial generators $X_d$ as being represented by 
\[
   \C P^d \longrightarrow \BU  
\]
classifying the hyperplane  line bundle over $\C P^d$. 
Altogether (using that $\MU_*$ commutes with direct limits), we get an isomorphism 
of graded $\MU_*$-algebras
\[
      \MU_*(B) \otimes A_*(G) \cong \MU_*[e_V, e_V^{-1}, Y_{V,d}]
\]
where $V$ runs over the nontrivial irreducible $G$-representations 
and $Y_{V,d}$, $|V|+1 \leq d < \infty$,   is the image of the element
\[
    X_{d-|V|} \otimes e_V^{-1} \in  \MU_{2d-2|V|}(\BU) \otimes A_{2|V|}(G)
\]
under the inclusion map $\BU \to B$ which embeds $\BU$ as the $V$th factor. 
If we had defined $X_0$ as the unit in $\MU_0(BU)$,
the corresponding element $Y_{V,|V|}$ would be equal to 
$e_V^{-1}$.

Heuristically, for a geometric class $[M] \in \Omega^G_*$, 
the element 
\[
   \phi_{\MU} \circ \Psi ([M]) \in \MU_*[e_V, e_V^{-1}, Y_{V,d}]
\]
describes the stable normal bundle of $M$ restricted to $M^G$. In 
particular, for each one-dimensional nontrivial irreducible complex $G$-representation $V$, we have
\[
 \iota \circ \phi_{\MU} \circ \Psi([{\rm P}(\C^d \oplus V)])  = 
 Y_{V,d} + e_{V^*}^{-d} \, ,
\]
where $V^*$ is the conjugate representation, $P$ denotes the projectivization and 
\[
   \iota : \MU_*[e_V, e_V^{-1}, Y_{V,d}]  \to  \MU_*[e_V, e_V^{-1}, Y_{V,d}] 
\]
is induced by the inverse of the $H$-space $B$ (we identify 
$\MU_*(B) \otimes A_*(G)$ and $\MU_*[e_V,e_V^{-1}, Y_{V,d}]$ 
by the isomorphism constructed before). The map $\iota$
interchanges the roles of the stable normal bundle of $M$ 
restricted to $M^G$ and of the normal bundle of $M^G$ in $M$.  Note
that $\iota \circ \Phi_{\MU}(e_V) = e_V$ for all $V$. 
The above equation  can be checked by 
inspecting the normal bundle of ${\rm P}(\C^d \oplus V)^G$ in 
${\rm P}(\C^d \oplus V)$, cf. Prop. 4.14. in \cite{S1} (in this 
reference, the map $\iota$ is mistakenly left out).    

We now discuss the geometric analogue of the map $\phi_{\MU}$
and define a map of $\MU_*$-algebras 
\[
  \phi_{\Omega}: \Omega_*^G \to  \MU_*(B) \otimes A_*(G) 
\]
in the following way: Let $M^n$ be a stable almost complex $G$-manifold and 
let 
\[
    F \subset M^G 
\]
be a connected component of the fixed point set. The normal 
bundle $\nu(F)$ of $F$ in $M$ is a complex $G$-bundle in a canonical 
way. Let $k$ be its 
complex dimension. Write 
\[
     \nu(F) = (E_{1} \otimes V_{1}) \oplus \ldots \oplus (E_{j} \otimes V_{j})
\]
with complex vector bundles 
\[
   E_1, \ldots, E_j
\]
and irreducible $G$-representations 
\[
   V_1, \ldots, V_j \, . 
\] 
Now define  
\[
  b_F := \overline{b}_F \otimes \big( e_{V_{1}}^{-|E_{1}|} \cdot \ldots \cdot e_{V_{j}}^{-|E_{j}|} 
  \big)\in \MU_{n - 2k }(B) \otimes A_{2k}(G)
\]
where $\overline{b}_F \in \MU_{n-2 k}(B)$ is represented by the map 
\[
   F \to B
\]
with the $V_i$-th component classifying the bundle $E_i$. Finally, we set 
\[
   \phi_{\Omega}([M]) := \sum_{F \subset M^G}  b_F \in (\MU(B) \otimes A(G))_n \, .
\]
The following result of tom Dieck shows that the geometric and 
homotopy theoretic fixed point maps are compatible with respect
to the Pontrjagin-Thom map. 

\begin{prop}[\cite{tD1}, Proposition 4.1] \label{tomD} The following diagram 
is commutative. 
\[
   \begin{CD} 
       \Omega^G_* @> \phi_{\Omega} >> \MU_*[e_V, e_V^{-1}, Y_{V,d}] \\
            @V\Psi VV                      @V {\rm id} VV             \\
        \MU^G_*    @> \iota \circ \phi_{\MU} >> \MU_*[e_V, e_V^{-1}, Y_{V,d}] 
   \end{CD}
\]
\end{prop}

\section{Geometric realizability} 

\begin{defn} We define the {\em geometric cone} in $\MU_*[e_V, e_V^{-1}, Y_{V,d}]$ 
as the $\MU_*$-subalgebra
\[   
    \Gamma_* := \MU_*[e_V^{-1}, Y_{V,d}] \, .
\]
\end{defn}

The following proposition justifies this name.   

\begin{prop} \label{simple} $\phi_{\Omega} (\Omega_*^G) \subset \Gamma_*$.
\end{prop}

\begin{proof} Due to the 
definition of the polynomial generators $Y_{V,d}$, the proof 
is nontrivial. Let $M$ be a
stable almost complex $G$-manifold and let 
\[
  F \subset M^G
\]
be a fixed point 
component with normal bundle $\nu(F)$ of complex dimension $k$. As above, we write 
\[
     \nu(F) =  (E_1 \otimes V_1) \oplus \ldots \oplus (E_j \otimes V_j) \, . 
\]
We show that the element   
\[
  \overline{b}_F \in \MU_{n-2k}(B) \otimes \Z 
\subset \MU_{n-2k}(B) \otimes A_{2k}(G) = \MU_*[e_{V}, e_V^{-1}, Y_{V,d}]
\]
is a polynomial in $e_{V_i}$ of degree at most  $| E_i |$ for all 
$V_i \in \{V_1, \ldots, V_j\}$. This in turn implies that 
\[
 b_F := \overline{b}_F \otimes \big( e_{V_{1}}^{-|E_{1}|} \cdot \ldots \cdot e_{V_{j}}^{-|E_{j}|} 
  \big)\in \MU_{n - 2k }(B) \otimes A_{2k}(G)
\]
is indeed an element of $\MU_*[e_V^{-1}, Y_{V,d}]$. 

The bordism class $\overline{b}_F$ can be regarded as being 
represented by a map   
\[
      F \ra \BU(|E_1|) \times \ldots \times \BU(|E_j|) \, . 
\]
Each element in $\MU_*(\BU(|E_i|))$, $1 \leq i \leq j$,  can be written 
as a sum of 
monomials in the $X_d$'s each of which contains at most $|E_i|$ 
factors of the form $X_d$, see e.g. \cite{Kochman}, Proposition 4.3.3 a). 
This implies that under the change of variables sending $X_d$ 
to $Y_{V_i,d+|V_i|} \otimes e_{V_i}$, each such element is 
expressible as a polynomial in $\MU_*[e_{V_i}, Y_{V_i,d}]$
of degree at most $|E_i|$ in $e_{V_i}$. Hence, using 
the   K\"unneth formula for computing 
\[
  \MU_*\Big( \BU(|E_1|) \times \ldots \times \BU(|E_j|) \Big) \, , 
\]
the element $\overline{b}_F$ is indeed of the required form. 
This finishes the proof of Proposition \ref{simple}. 
\end{proof}

\begin{thm} \label{main} Let $G = S^1 \times \ldots \times S^1$. 
Then the  commutative diagram of $\MU_*$-algebras
\[
   \begin{CD} 
           \Omega_*^G @> \phi_{\Omega} >> \Gamma_*     \\
                @V \Psi VV        @V {\rm incl.} VV \\
           \MU_*^G @> \iota \circ \phi_{\MU} >> \MU_*[e_V, e_V^{-1}, Y_{V,d}] 
   \end{CD}
\]
induced by the commutative fixed point square of tom Dieck (see 
Proposition \ref{tomD}) is a pull back square. All maps 
in this diagram are injective. 
\end{thm}

For $G$ a torus, the map $\phi_{\MU}$ is injective 
by \cite{S1}, Proposition 4.5 and Corollary 5.2. Thus (together 
with Proposition \ref{loffler}) only the 
pull back property requires proof. 

The last theorem implies the following 
relation of geometric and homotopy theoretic equivariant 
bordism. 

\begin{cor} \label{optimal} Let $G$ be a torus. Then the Pontrjagin-Thom 
map induces an isomorphism of $\MU_*$-algebras
\[
     \Omega_*^G \cong \im(\iota \circ \phi_{\MU}) \cap \Gamma_* \, .
\]
\end{cor}

Hence, the coefficients of the geometric equivariant $(S^1)^r$-bordism ring 
constitute a certain subring of $\MU_*[e_V, e_V^{-1}, Y_{V,d}]$ whose 
elements are characterized by the following two conditions: 
Firstly, they must lie in the image of $\iota \circ \phi_{\MU}$. Heuristically,
the local fixed point data specified 
by an element in $\MU_*[e_V, e_V^{-1}, Y_{V,d}]$ can be ``closed up''  (without 
introducing new fixed points) such as to define  
an actual class in $\MU^G_*$. This property can be formalized
using bordism with respect to families of subgroups of $G$ (see below). 
Secondly, they must lie in the geometric cone.

The first of these two conditions has been extensively studied by tom Dieck 
\cite{tD1,tD3} in terms of integrality conditions 
related to the localization techniques of  Atiyah-Segal.
We will come back to this description in Section \ref{integral}.

We now turn to the proof of Theorem \ref{main}. 
Recall that a {\em family of subgroups} $\mathcal{F}$ of a topological 
group $G$ is a set of closed subgroups of  $G$ which is closed under
conjugation and under taking subgroups. Such 
a family is supposed to encode the possible isotropy groups occuring 
in a given $G$-space. We define special families of subgroups of $G$:
\begin{eqnarray*}
     \mathcal{A} & := & \{ H < G \} \, , \\
     \mathcal{P} & := & \{ H < G ~|~ H \neq G \} \, .
\end{eqnarray*}
From now on, we assume that $G$ is a compact Lie group. For each family of 
subgroups $\mathcal{F}$, there exists a classifying space $E \mathcal{F}$ 
\cite{tD4}, I.(6.6), which is a terminal object in the homotopy category 
of $\mathcal{F}$-numerable  $G$-spaces. The space $E \mathcal{F}$ 
is characterized by  the properties \begin{itemize}
   \item[-] $E \mathcal{F}^H \simeq *$, if $H \in \mathcal{F}$\, ,  
   \item[-] $E \mathcal{F}^H = \emptyset$, if $H \notin \mathcal{F}$. 
\end{itemize}
(using nonequivariant homotopy equivalences). Furthermore, it is unique up to $G$-homotopy equivalence. For a $G$-space $X$ and 
a pair of families of subgroups $(\mathcal{F}, \mathcal{F}')$ (i.e.
$\mathcal{F}' \subset \mathcal{F}$), we set 
\begin{eqnarray*}
    \Omega_n^G[\mathcal{F},\mathcal{F}'](X) & := & \Omega_n^G( X \times E \mathcal{F}, X \times E\mathcal{F}')\, ,  \\
   \MU^G_n[\mathcal{F}, \mathcal{F}'](X) & : = & \MU_n^G(X \times E
 \mathcal{F}, X \times E \mathcal{F}')\, .
\end{eqnarray*}
cf. \cite{C}, p. 339. The groups $\Omega^G_n[\mathcal{F}, \mathcal{F}'](X)$ 
consist of $G$-bordism classes  of stable almost complex $G$-manifolds $(M^n, \partial M)$ with 
boundary and with reference maps to $X$ such that 
all isotropy groups occuring in $M$ (resp. in $\partial M$) 
lie in $\mathcal{F}$ (resp. in $\mathcal{F}'$). 
The long exact sequence of the pair 
\[
   (X \times E \mathcal{F}, X \times E \mathcal{F}')
\] 
is the Conner-Floyd exact sequence
\[
    \ldots \to \Omega^G_{n}[\mathcal{F}'](X) \to \Omega^G_n[\mathcal{F}] (X) \to 
    \Omega_{n}^G[\mathcal{F}, \mathcal{F}'] (X) \to 
   \Omega^G_{n-1}[\mathcal{F}'](X) \to \ldots
\]
The proof of Theorem \ref{main} is based on studying the commutative diagram with 
exact rows
\[
  \begin{CD} 
      \ldots @>>>   \Omega^G_n @>>> \Omega^G_n[\mathcal{A}, \mathcal{P}] @>>> \Omega^G_{n-1}[\mathcal{P}] @>>> \ldots \\
         @VVV          @V\Psi = \Psi_{[\mathcal{A}]}VV                @V\Psi_{[\mathcal{A}, \mathcal{P}]} VV        @V\Psi_{[\mathcal{P}]}VV                 @VVV \\
      \ldots @>>>   \MU^G_n    @>>> \MU^G_n[\mathcal{A}, \mathcal{P}]   @>>> \MU^G_{n-1}[\mathcal{P}] @>>> \ldots
  \end{CD}
\]
whose vertical arrows are Pontrjagin-Thom maps. We know from \cite{C}, 
Theorem 5.4, that $\Psi_{[\mathcal{A}]}$ is injective. In a first step, we will identify the 
terms occuring in the second column by showing the existence of a commutative diagram 
\[
    \begin{CD} 
        \Omega^G_*[\mathcal{A}, \mathcal{P}] @> \cong >> \Gamma_* \\
      @V\Psi_{[\mathcal{A}, \mathcal{P}]}VV             @V {\rm incl.} VV   \\
        \MU^G_*[\mathcal{A}, \mathcal{P}] @> \cong >> \MU_*[e_V, e_{V}^{-1}, Y_{V,d}] 
    \end{CD}
\]
such that the compositions of the horizontal maps with the canonical maps 
\[
  \Omega^G_n  \to  \Omega^G_n[\mathcal{A}, \mathcal{P}] \, , ~~~
    \MU^G_n     \to  \MU^G_n[\mathcal{A}, \mathcal{P}] 
\]
coincide with the fixed point maps $\phi_{\Omega}$ and 
$\iota \circ \phi_{\MU}$.  
This will be done 
in Proposition \ref{later}. In 
a second step we show in Proposition \ref{PT} that the map 
$\Psi_{[\mathcal{P}]}$ is injective.  For this, we 
use equivariant characteristic numbers as introduced by tom Dieck
\cite{tD1a,tD2}.
From these facts together with the injectivity of $\Psi_{[\mathcal{A}]}$ and 
the injectivity of the fixed point map $\iota \circ \phi_{\MU}$, the assertion 
of Theorem \ref{main} follows by an easy diagram chase. 

\begin{prop} \label{later} Let $G$ be a compact Lie group. Then 
there are  ring isomorphisms 
\[
   \kappa_{\Omega} : \Omega^G_*[\mathcal{A}, \mathcal{P}] \to \Gamma_* 
\]
and 
\[
   \kappa_{\MU}: \MU^G_*[\mathcal{A}, \mathcal{P}] \to \MU_*[e_V, e_V^{-1}, Y_{V,d}]
\]
such that  the diagram 
\[
 \begin{CD} 
 \Omega^G_*[\mathcal{A}, \mathcal{P}] @> \kappa_{\Omega} >>  \Gamma_*  \\ 
                  @V\Psi_{[\mathcal{A}, \mathcal{P}]}VV        @V{\rm incl.} VV \\
 \MU_{*}^G[\mathcal{A}, \mathcal{P}]     @> \iota \circ \kappa_{\MU}  >> \MU_*[e_V, e_V^{-1}, Y_{V,d}] \\
   \end{CD}
\]
is commutative and the compositions
\[
 \Omega_{*}^G \to \Omega^G_*[\mathcal{A}, \mathcal{P}] \stackrel{\kappa_{\Omega}}{\longrightarrow}
 \Gamma_*
\]
and 
\[
    \MU^G_* \to \MU^G_*[\mathcal{A}, \mathcal{P}] \stackrel{\kappa_{\MU}}{\longrightarrow} 
\MU^G_*[e_V, e_V^{-1}, Y_{V,d}]
\]
coincide with the fixed point maps $\phi_{\Omega}$ and $\phi_{\MU}$, respectively. 
\end{prop}

\begin{proof} We construct a map 
\[
   \kappa_{\Omega} : \Omega^G_*[\mathcal{A}, \mathcal{P}] \to \Gamma_* \subset \MU_*(B) \otimes A_*(G)  
\]
as follows. The group $\Omega^G_n[\mathcal{A}, \mathcal{P}]$ is generated by
elements represented by  (disc bundles of) unitary $G$-bundles 
\[
    E \to M 
\]
where $M$ is a stable almost complex manifold of dimension $n - 2k$ 
with trivial $G$-action, $E$ has complex fibre dimension $k$ and 
$E^G = M$. We 
define the image of $[E \to M]$ under $\kappa_{\Omega}$
in exactly the same way as 
we defined the geometric fixed point map $\phi_{\Omega}$ above. The proof of 
the fact that we get indeed an element in $\Gamma_*$ is similar as before. 
An inverse of $\kappa_{\Omega}$ is constructed as follows. Let $V$
be a nontrivial irreducible unitary $G$-representation. We map 
$Y_{V,d}$ to the disc bundle of the $G$-bundle 
\[   
    E \otimes V \to \C P^{d-|V|}
\]
where $E \to \C P^{d-|V|}$ is the hyperplane line bundle. 
The element $e_{V}^{-1}$ 
is mapped to the disc bundle of the trivial $G$-bundle 
\[
   V \to *  \, . 
\]
Now  $\kappa_{\Omega}^{-1}$ is determined by the 
fact that it is a map of $\MU_*$-algebras.

On the homotopy theoretic side, the construction is as follows. 
The map 
\[
  \kappa_{\MU} :  \MU_{*}^G[\mathcal{A}, \mathcal{P}] \to \MU_*[e_V, e_V^{-1}, Y_{V,d}]
\]
is defined as the composition  
\[
  \MU_*^G[\mathcal{A}, \mathcal{P}] \cong \widetilde{\MU}_*^G(\Sigma E \mathcal{P})
 \cong (\phi^G \MU^G)_* \cong \MU_*[e_V, e_V^{-1}, Y_{V,d}] \, .
\]
In this sequence, the first isomorphism is induced by the excisive
inclusion of $G$-spaces
\[
 (E \mathcal{A}, E\mathcal{P}) \simeq 
(C_- E\mathcal{P} , E \mathcal{P}) \to (\Sigma E \mathcal{P}, 
C_+ E \mathcal{P})
\, , 
\]
where $C_-$ and $C_+$ denote the lower and upper cone in the unreduced 
suspension 
\[
  \Sigma E \mathcal{P}= ([0,1] \times E \mathcal{P}) / \sim
\]
which we consider as being equipped 
with the basepoint $[\{1\} \times E \mathcal{P}]$ and as containing 
$E\mathcal{P}$ as the subspace $\{1/2\} \times E \mathcal{P}$. The second 
isomorphism is induced by Lemma 4.2. in \cite{S1} with $Z := \Sigma E \mathcal{P}$
and the third isomorphism is Theorem 4.10. in {\em loc.~cit.} 

The commutativity of the diagram in Proposition \ref{later} is proven 
in a similar fashion as in the case of tom Dieck's fixed point square
(cf.~Proposition \ref{tomD}).
We omit the details.  The remaining assertions are immediate. 
\end{proof}

We now deal with the second step of the proof of Theorem \ref{main} and show 
the following variant of \cite{C}, Theorem 5.4.

\begin{prop} \label{PT} Let $G$ be topologically cyclic group, i.e. $G$ 
is of the form $S^1 \times \ldots \times S^1 \times \Z/k$, 
$k = 1, 2, \ldots$. Then the 
Pontrjagin-Thom map
\[
    \Psi_{[\mathcal{P}]} : \Omega_n^G[\mathcal{P}] \to \MU_n^G[\mathcal{P}]
\]
is injective for all $n$. 
\end{prop}

\begin{proof} If $G = S^1$, the assertion is implied 
by  \cite{C}, 
Theorem 7.1.(3) with the choice $G  := \{e\}$. The general case 
can conveniently be dealt with using equivariant 
characteristic numbers \cite{tD1a,tD2}. It follows from \cite{tD2}, 
Theorem 3, that 
for topologically cyclic $G$ and a unitary  $G$-representation $W$, the characteristic number map 
\[
    \Omega^G_{*}(SW) \to K_G^{-*}(SW)[[a_1, a_2, \ldots ]]
\]
is injective. Here, $SW$ denotes the unit sphere in $W$. But this map 
factors through the Pontrjagin-Thom map 
\[
    \Omega^G_{*}(SW) \to \MU^G_*(SW) \, . 
\]
Now observe that  
\[ 
   E \mathcal{P} = \lim_{\longrightarrow_W} SW
\]
where the colimit is taken with respect to a directed set of complex $G$-modules 
without trivial direct summands and so that ultimately each nontrivial irreducible appears with arbitrarily large multiplicity in some 
$G$-module $W$. Finally, the claim of Proposition \ref{PT} follows
because  homology commutes with direct limits: 
\begin{eqnarray*} 
    \Omega^G_*(E \mathcal{P}) & = & \lim_{\longrightarrow_W} \Omega^G_{*}(SW)  \, , \\ 
    \MU^G_*(E \mathcal{P}) & = & \lim_{\longrightarrow_W} \Omega^G_{*}(SW) \, .
\end{eqnarray*}
\end{proof} 

We remark that Proposition \ref{PT} holds for any compact abelian Lie group $G$
and with $\mathcal{P}$ replaced by any family of subgroups of $G$. 
The proof of this general statement can be carried out using 
similar ideas as in \cite{C} 
and is based on a systematic use of the Conner-Floyd exact sequences
for bordism with respect to  families of subgroups. However, it 
is considerably more tedious than the proof presented above for 
topologically cyclic $G$. Since only this case is 
of importance for us, we decided not to include
the more general case in this paper.

\section{The case of nonconnected $G$}

There are simple  examples showing that if $G$ is not connected
and not of prime order, then 
the geometric 
realizability of an element in $\MU^{G}_*$ is not determined 
by the normal 
data around fixed point sets. Let $n$ be a 
natural number, $G = \Z / n \times \Z/n$ and let 
\[
   f : \Z/n \to \BU(1,G)
\]
classify the bundle 
\[
     \Z/n \times V \to \Z/n
\]
where $V$ is the canonical complex $1$-dimensional representation of $\Z/n$, 
the space $\Z/n \times V$ carries 
the product $(\Z/n \times \Z/n)$-action and the bundle map 
is simply projection onto the first factor. 
The map
\[
    \Z/n \times S^0 \subset \Z/n \times S^V \to T(\xi^{G}_{1})
\]
induced by $f$ represents an element $c \in \MU_{-2}^G$ which is different from zero, because 
its restriction to $(\{1\} \times \Z/n)$-bordism is $n$ times the 
Euler class $e_V \in \MU^{\Z/n}_{-2}$ (this is known to 
be a nontorsion class). However,  
\[
   \phi_{\MU}(c) = 0 \, . 
\]
A similar example exists for $G = \Z / n^2$ with $f:\Z/n \to \BU(1, G)$ 
classifying the bundle 
\[
    \Z/n \times W \to \Z/n
\]
where $W$ is the canonical one dimensional $\Z/n^2$-representation and $G$ acts
on the $\Z/n$-factors via the projection $G \to \Z/n$ with kernel 
$\Z/n$.  
The restriction of the corresponding bordism class  to 
$\Z/n$-bordism 
(where $\Z/n \subset G$ is of index $n$), is 
again given as $n$ times the Euler class $e_V \in \MU^{\Z/n}_{-2}$. 

One reason of the failure of Theorem \ref{main} for arbitrary 
compact abelian $G$ is the fact
that the fixed point map 
\[
  \phi_{\MU} : \MU^G_* \to \MU_*[e_V, e_V^{-1}, Y_{V,d}]
\]
is injective if and only if $G$ is a torus, see \cite{S1}, Theorem 5.1. 
It seems plausible that for geometric 
realizability of an element in $\MU_*^G$ for general $G$, a 
series of obstructions must vanish each one of which lies in some 
group 
\[
   \MU^G_*[\mathcal{F}, \mathcal{F}']
\]
with an adjacent pair of families $(\mathcal{F}, \mathcal{F}')$ and each 
one of which is only defined if the previous one vanishes. However, 
we do not know a concise statement in this direction.

\section{Integrality} \label{integral}

We will study the image of the map  
\[
    \MU^G _* \to \MU_*^G[\mathcal{A}, \mathcal{P}] 
   \stackrel{\iota \circ \kappa_{\MU}}\cong \MU_*[e_V, e_V^{-1}, Y_{V,d}] 
\]
(which occurs in the lower line of the diagram in Theorem \ref{main})
in the light of integrality conditions formulated by tom Dieck. We have a natural 
multiplicative transformation  
\[
    \eta : \MU_G^*(X) \to \MU^*(EG \times_G X) 
\]
where $X$ is an arbitrary $G$-CW complex. This map (called ``bundling transformation'' in the papers of tom Dieck) is induced by the 
projection 
\[
    EG \times X \to X
\]
yielding a map 
\[
   \MU^*_G(X) \to \MU^*_G(EG \times X) \cong \MU^*( EG \times_G X) \, . 
\]
For a complex $G$-representation $W$, let $e(W)\in \MU^2(BG)$ denote 
the bordism theoretic Euler class of the complex vector bundle 
\[
    EG \times_G W \to BG \, . 
\]
It follows directly from the definition of $\eta$ that 
\[
  \eta(e_W) = e(W) \, . 
\]
Here, we regard $e_W$ as an element in $\MU^{2 |W|}_G = \MU_{-2 |W|}^G$. 
Recall that there is an isomorphism 
\[
   \MU^*(B(S^1)^r ) \cong \MU^*[[C_1, \ldots, C_r]]
\]
where $C_i \in \MU^2(B(S^1)^r)$ is the Euler class of the 
$(S^1)^r$-representation 
$\rho_i$ which is induced from the standard one-dimensional 
$S^1$-representation by projecting 
$(S^1)^r$ onto the $i$th factor.   If 
\[
  V = \rho_1^{\otimes \mu_1} \otimes \ldots \otimes \rho_r^{\otimes \mu_r}\, , 
\]
then 
\[
   e(V) = [\mu_1]_F C_1 +_F \ldots +_F [\mu_r]_F C_r
\]
with Lazard's universal formal group law $+_F$ in $\MU^*[[C_1, \ldots, C_r]]$. 

Let $S \subset \MU_G^* = \MU^G_{-*}$ denote the multiplicative subset of Euler classes of 
$G$-representations without trivial direct summand. We denote 
the corresponding subset of $\MU^*(BG)$ by $S$ as well and obtain an induced
map 
\[
 S^{-1} \eta : \MU^G_*[\mathcal{A}, \mathcal{P}] \cong S^{-1} \MU^G_* \to S^{-1} \MU^{-*}(BG) 
\]
compatible with $\eta$ in the obvious way. 
The following result  of tom Dieck \cite{tD2}  characterizes 
the image of 
\[
  \iota \circ \Phi_{\MU} : \MU^G_* \to \MU^G_*[\mathcal{A}, \mathcal{P}] 
   \stackrel{\iota \circ \kappa_{\MU}}{\cong} \MU_*[e_V, e_V^{-1}, Y_{V,d}] \, . 
\] 
We will regard $S^{-1} \eta$ as 
being defined on $\MU_*[e_V, e_V^{-1}, Y_{V,d}]$ by identifying this 
module with $\MU^G_*[\mathcal{A}, \mathcal{P}]$ using the isomorphism 
$\iota \circ \kappa_{\MU}$.  

\begin{prop}[\cite{tD2}] Let $G$ be a topologically cyclic group and let 
\[
    \lambda : \MU^*(BG) \to S^{-1} \MU^*(BG) 
\]
be the localization map. Then 
\[
   \im (\iota \circ \Phi_{\MU})  = \{ x \in \MU_*[e_V, e_V^{-1}, Y_{V,d}] ~|~ 
S^{-1} \eta(x) \in \lambda ( \MU^{-*}(BG)) \}
\]
\end{prop} 

For $G = S^1 \times \ldots \times S^1$, the localization map 
$\lambda$ is injective, because $\MU^*(BG)$ has 
no zero divisors. Therefore, $\MU^*[[C_1, \ldots , C_r]]$ 
can be considered as a subring of its $S$-localization. 
In combination with Theorem \ref{main} and by the fact that 
for toral $G$ the map $\eta$ is injective (because the characteristic 
number map factors through $\eta$), we now get the following
description of the geometric unitary $G$-bordism ring. 

\begin{thm} Let $G = (S^1)^r$.  Let 
\[
   {\rm incl.} : \MU_{*}[e_{V}^{-1}, Y_{V,d}] \to \MU_{*}[e_V, e_{V}^{-1}, Y_{V,d}] 
\]
be the canonical inclusion. As before, we identify $\MU_n$ and $\MU^{-n}$. Then 
the diagram  
\[
  \begin{CD} 
        \Omega^{G}_* @>\Phi_{\Omega}>> \MU^{-*}[e_{V}^{-1}, Y_{V,d}] \\
             @V\eta \circ \Psi VV           @V (S^{-1} \eta) \circ \iota \circ {\rm incl.}VV \\
         \MU^{-*}[[C_1, \ldots, C_r]] @>\lambda >>  S^{-1} \MU^{-*}[[C_1, \ldots, C_r]]
  \end{CD}
\]
is a pull back square of $\MU_*$-algebras. 
\end{thm}

Acknowledgements: The author is grateful to U. Bunke and T. Schick 
for helpful comments.

\end{document}